\documentclass[a4paper]{article}
\usepackage[latin1]{inputenc}
\usepackage[english]{babel}
\usepackage{amssymb,amsmath,amscd,amstext}
\usepackage{amsfonts}
\usepackage{amsthm}
\usepackage[dvips]{graphicx}
\usepackage{epsfig}
\usepackage{palatino}
\usepackage{xypic}
\input{xy}
\xyoption{all}
\renewcommand{\vec}[1]{\boldsymbol{#1}}

\newcounter{eno}
\newcounter{ch}
\newcounter{sec}
\newcounter{no}
\def \sec#1 {\stepcounter{sec} \bigskip \par{\noindent \begin{center}{\bf \S \arabic{sec} {\bf #1}  \setcounter{eno}{1} \setcounter{no}{1}} \end{center}}}
\def \prop#1 {\bigskip\par{\noindent \bf Proposition \arabic{sec}.\arabic{no}} \stepcounter{no} { \it #1}}
\def \cor#1 {\bigskip\par{\noindent \bf Corollary \arabic{sec}.\arabic{no}} \stepcounter{no} {\it #1}}
\def \lem#1 {\bigskip\par{\noindent  \bf Lemma \arabic{sec}.\arabic{no}} \stepcounter{no} { \it #1}}
\def \thm#1 {\bigskip\par{\noindent \bf Theorem \arabic{sec}.\arabic{no}} \stepcounter{no} { \it #1}}
\def \conj#1 {\bigskip\par{\noindent \bf Conjecture \arabic{sec}.\arabic{no}} \stepcounter{no} { \it #1}}
\def \define#1 {\bigskip\par{\noindent \bf Definition \arabic{sec}.\arabic{no}} \stepcounter{no} { \it #1} \hspace{2mm}}
\def \egg#1 {\bigskip\par{\noindent \bf Example \arabic{sec}.\arabic{no}} \stepcounter{no} {\it #1}}

\newcommand {\R} {\mathbb R}
\newcommand {\C} {\mathbb C}
\newcommand {\Z} {\mathbb Z}
\newcommand {\Q} {\mathbb Q}
\newcommand {\ZG} {\mathbb {Z}[G]}
\newcommand {\QG} {\mathbb {Q}[G]}

\def \proof#1 { Proof:$\,\,\,$  #1  ${}$\hfill $\Box$
\,\, \newline}

\def\cs#1{\texttt{\char`\\#1}}

\def\describecs{%
  \ifvmode\medskip\noindent\fi
  \hspace*{-1cm}\cs}
\begin{document}{\rm

\title{The third homotopy module of a 2-complex}
\date{}
\author{W.H.Mannan}
\maketitle

\catcode`\|\active
\def|#1|{\textnormal{\small\itshape#1}}

\vspace{-5mm}

\begin{quote}{\small \noindent {\bf A{\tiny BSTRACT}:}  Given a connected 2- complex
$X$ with fundamental group $G$, we show how $\pi_3(X)$ may be
computed as a module over $\ZG$. Further we show that if $X$ is a
finite connected 2- complex with $G\,(=\pi_1(X))$ finite of odd
order then the stable class of $\pi_3(X)$ is determined by $G$. In
this case we also show that $\pi_3(X) \otimes \Q$ is free over
$\QG$.}\end{quote}

{\bf MSC}: 55Q15, 55Q20 55Q91 \hfill {\bf Keywords} Homotopy groups, 2-complexes

\bigskip
\noindent {\bf email}: \verb|wajid@mannan.info|

\sec{Introduction}

Let $X$ be a connected 2- complex with fundamental group $G$. Let
$C_*(\tilde{X})$  denote the algebraic complex over $\ZG$
associated to the universal cover of $X$.  We may write this
algebraic complex:
$$
F_2 \stackrel{\partial_2}{\to} F_1 \stackrel{\partial_1}{\to} F_0
$$

\noindent From the Hurewicz isomorphism theorem, we know that
$$\pi_2(X) \cong H_2(C_*(\tilde{X});\ZG) \cong {\rm
Ker}\partial_2$$ as (right) modules over $\ZG$.  Hence $\pi_2(X)$
is easily attained from $C_*(\tilde{X})$.

\bigskip
Now let $J$ be a module over $\ZG$, which is freely generated over
$\Z$.  There is a $G$- action on $J \otimes_{\Z} J$, given by $(a
\otimes b)g = ag \otimes bg$, making it a $\ZG$ module. Let
$\sigma$ be the $\ZG$- linear automorphism of $J \otimes_{\Z} J$
defined by $\sigma(a \otimes b) =b \otimes a$.

\define{} $S^2(J)=\lbrace x \in J \otimes J \vert \,\, \sigma x=x \rbrace $

\bigskip
We now state our main result:

\bigskip
\noindent{\bf Theorem A}  We have an isomorphism of $\ZG$ modules:
$\pi_3(X)\cong S^2(\pi_2(X))$

\bigskip
When $X$ is a finite connected 2- complex, Schanuel's lemma tells
us that the stable class of $\pi_2(X)$ is determined by $G$.  We
prove an analogous result for $\pi_3(X)$, under the added
assumption that $G$ is finite, of odd order:

\bigskip
\noindent {\bf Theorem B} {\it Let $X$ and $X'$ be finite
connected 2- complexes, with finite fundamental group $G$.  If the
order of $G$ is odd, then $\pi_3(X)$ and $\pi_3(X')$ are stably
equivalent.}

\bigskip
 Also if $G$ is finite, of odd order we have the following result:

\bigskip
\noindent {\bf Theorem C} {\it Let $X$ be a finite connected $2$-
complex with finite fundamental group $G$.  If the order of $G$ is
odd, then $\pi_3(X) \otimes \Q$ is free over $\QG$.}

\define{$G$- quadratic}  A map
$f:A \to B$ between $\ZG$ modules $A$, $B$, is $G$- quadratic if
given $a,b,c \in A$ and $g \in G$ we have:

i)\,\,$f(-a)=f(a)$

ii)\,\,$f(a+b+c)=f(a+b)+f(a+c)+f(b+c)-f(a)-f(b)-f(c)$

iii)\,\,$f(ag)=f(a)g$

\lem{Let $f:A \to B$ be $G$- quadratic;}

\hspace{2mm}i)\,\,$f(0)=0$

\hspace{2mm}ii)\,\, $f(a+b)+f(a-b)=2f(a)+2f(b)$

\bigskip
\proof{i)\,\,$f(0)=f(0+0+0)=0$
\newline
${}$\hspace{17mm}ii)\,\,$f(a+b)+f(a-b)=f(a+b)+f(a-b)+f(b-b)$\newline
${}$\hspace{21mm}$=f(a+b-b)+f(a)+f(b)+f(-b) =2f(a)+2f(b)$}

Take $J$ as before.  We have a $G$- quadratic map $q:J \to S^2(J)$
which maps $\alpha \mapsto \alpha \otimes \alpha$, for any $\alpha
\in J$.  The main argument of this paper will depend on a certain
universal property of $q$, which we mention as a preliminary.

\bigskip
Let $\{\alpha_i \vert\,\, i \in I\}$ be a basis for $J$ over $\Z$
for some set $I$.  Given $\alpha \in J$ we may write $\alpha =
\sum_{i\in I} \alpha_i \lambda_i$ with each $\lambda_i \in \Z$ and
only finitely many $\lambda_i$ non-zero. We define the norm of
$\alpha$ to be $\sum_{i\in I} \vert\lambda_i \vert$.

The $\{\alpha_i \otimes \alpha_j \vert \, i,j \in I\}$ are a basis
for $J \otimes_{\Z}J$ over $\Z$.  So given $\omega \in S^2(J)$, we
may write $\omega = \sum_{i,j \in I} \lambda_{ij} \alpha_i \otimes
\alpha_j$.  Clearly for each $i,j \in I$ we must have
$\lambda_{ij}=\lambda_{ji}$.  Hence $S^2(J)$ is generated freely
over $\Z$ by the $\alpha_i \otimes \alpha_i$ and the $\alpha_i
\otimes\alpha_j+\alpha_j \otimes \alpha_i$, $i \neq j$.

For each $i,j \in I$ we have $\alpha_i \otimes\alpha_j+\alpha_j
\otimes \alpha_i= q(\alpha_i +\alpha_j)-q(\alpha_i)+q(\alpha_j)$.
Hence another basis for $S^2(J)$ over $\Z$ is given by the
$q(\alpha_i)$ and the $q(\alpha_i+\alpha_j)$, $i \neq j$.

\lem{The map $q:J \to S^2(J)$ is universal: Given any $G$-
quadratic map $f:J \to M$, there exists a unique $\ZG$- linear map
$\phi:S^2(J) \to M$ making the following diagram commute:}

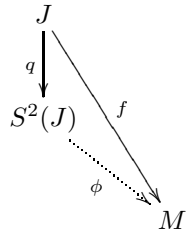
\begin{center}\center{\xymatrix{&&&&J\ar[d]_q
\ar[ddr]^f&\\&&&&S^2(J)\ar@{.>}[dr]_\phi\\&&&&& M}}
\end{center}

\proof{ We define $\phi$ linearly on a basis by:
\newline
${}$\hspace{71mm}$\phi q(\alpha_i)=f(\alpha_i)$
\newline
${}$\hspace{71mm}$\phi
q(\alpha_i+\alpha_j)=f(\alpha_i+\alpha_j),\,\,i \neq j$
\newline
Clearly no other map can make the diagram commute.  It remains to
show that $\phi q=f$ and that $\phi$ is $\ZG$- linear.
\newline
${}$\hspace{2mm}By lemma 1.3 ii) we have $\phi
q(\alpha_i-\alpha_j)= f(\alpha_i-\alpha_j)$.  Hence, by
construction we have $\phi q(\alpha)=f(\alpha)$ for any $\alpha
\in J$ with norm less than or equal to $2$.
\newline
${}$\hspace{2mm} We proceed by induction.  Suppose we have
verified that $\phi q(\alpha)=f(\alpha)$ for all $\alpha \in J$
with norm less than $r$, where $r \geq 3$. Then consider any
$\beta \in J$ with norm equal to $r$. We may write
$\beta=\alpha+x+y$ with $x=\pm\alpha_i$ and $y=\pm \alpha_j$ for
some $i, j \in I$, such that the norms of $\alpha$, $\alpha+x$,
$\alpha+y$ are all less than $r$.
\newline
${}$\hspace{2mm} Then property ii) of a $G$- quadratic map implies
that $$\phi q(\beta)=\phi q(\alpha+x+y)=f(\alpha+x+y)= f(\beta)$$
So we have verified that $\phi q(\alpha)=f(\alpha)$ for all
$\alpha \in J$ with norm less than $r+1$.  Hence by induction
$\phi q =f$.  It remains to show that $\phi$ is $\ZG$- linear.
\newline
${}$\hspace{2mm} By construction $\phi$ is $\Z$- linear.  Given
$\alpha \in J$ we have $\phi(q(\alpha)g) = \phi q(\alpha
g)=f(\alpha g)=f(\alpha)g=(\phi q(\alpha))g$.  Any element of
$S^2(J)$ is a $\Z$-linear combination of elements of the form
$q(\alpha)$, $\alpha \in J$. Hence $\phi$ is $\ZG$- linear on the
whole of $S^2(J)$.}

In \S2 we consider a general topological space $X$ with $\pi_1(X)
= G$. We give a construction for a $G$- quadratic map $\pi_2(X)
\to \pi_3(X)$.  Whenever $\pi_2(X)$ is freely generated over $\Z$
lemma 1.4 will therefore give us a $\ZG$- linear map
$\phi:S^2(\pi_2(X)) \to \pi_3(X)$.  In \S3 we show that when $X$
is a connected 2- complex, this map is an isomorphism.

\sec{The Hopf fibration}

Let $X$ be an arbitrary topological space with basepoint $*$. Let
$G$ denote $\pi_1(X)$.  We will consider the map which the Hopf
fibration $h:S^3 \to S^2$ induces on homotopy: $h^*:\pi_2(X) \to
\pi_3(X)$.  In particular we will show that it is $G$- quadratic.

\bigskip
We may regard $S^3$ as $\{(z,w) \in \C ^2 \vert\,z \bar{z}+w
\bar{w}=1\}$. Define a relation $\sim$, by setting $(z,w) \sim (z
\lambda, w\lambda)$ whenever $\lambda \in \C$ satisfies $\lambda
\bar{\lambda} =1$.

\define{Hopf fibration} The Hopf fibration, $h:S^3 \to S^2$ is the natural map
$S^3 \to S^3/\sim$ composed with the identifications: $S^3 /\sim$
$\cong \C P^1 \cong S^2$.

\bigskip
We take $(0,1)$ as basepoint of $S^3$ and let ${\tiny +}= h(0,1)
\in S^2$ be the basepoint of $S^2$.  We define $h^*:\pi_2(X) \to
\pi_3(X)$ to be the map given by precomposition with $h$. As
composition of maps is associative we know that $h^*$ is natural:

\lem{Let $f:X \to Y$ be a basepoint preserving continuous map of
based topological spaces.  Then $h^*f_*=f_*h^*$.}

\bigskip
We first show that $h^*$ satisfies property i) of a $G$- quadratic
map:

\lem{Let $\alpha \in \pi_2(X)$.  Then $h^*(-\alpha) =
h^*(\alpha)$.}

\proof{Let $R:S^3 \to S^3$ be the rotation given by $(z,w) \mapsto
(\bar{z},\bar{w})$.  Let $\hat{R}:S^2 \to S^2$ be the reflection
induced on $S^3 / \sim$.  Then $R$ is homotopic to the identity
(keeping the basepoint $(0,1)$ fixed).  So as elements of
$\pi_3(X)$ we have:$$ h^*(-\alpha)=\alpha \circ \hat{R} \circ h
=\alpha \circ h \circ R = \alpha \circ h =h^*(\alpha)$$}

Given $\alpha, \beta \in \pi_2(X)$, we follow \cite{Whit}, Chap.
X, \S7 by denoting their Whitehead product $[\alpha, \beta]$.  We
also recall the following results from that section:

\prop{The Whitehead product satisfies:
\newline
i) (Naturality, see (7.2)) Given $\alpha, \beta \in \pi_2(X)$ we
have $f_*[\alpha,\beta]=[f_*\alpha,f_*\beta]$ for any continuous
basepoint preserving map $f:X \to Y$.
\newline
ii) (Commutativity, see (7.5)) Given $\alpha, \beta \in \pi_2(X)$
we have $[\alpha, \beta]= [\beta,\alpha]$.
\newline
iii) (Biadditivity, see (7.12)) Given $\alpha_1,\alpha_2, \beta
\in \pi_2(X)$ we have: \newline
 $[\alpha_1+\alpha_2, \beta]
=[\alpha_1, \beta] + [\alpha_2, \beta]$ and $[\beta,
\alpha_1+\alpha_2] =[\beta, \alpha_1] + [\beta, \alpha_2]$.}

\bigskip
Whitehead products are related to $h^*$ by the following result:

\lem{$h^*(\alpha+\beta)=h^*(\alpha)+h^*(\beta)+[\alpha,\beta]$}

\proof{See \cite{Whit}, Chap. XI, (1.16) and from (4.4) note that
we may choose conventions so that the Hopf invariant of $h$ is
$1$.}

We may now prove that $h^*$ satisfies property ii) of a $G$-
quadratic map:

\lem{$h^*(\alpha+\beta+\gamma)=h^*(\alpha+\beta)+
h^*(\alpha+\gamma)+h^*(\beta+\gamma)$
\newline ${}$ \hspace{61mm} $-h^*(\alpha)-h^*(\beta)-h^*(\gamma)$}

\bigskip
\proof{$h^*(\alpha+\beta+\gamma)= h^*(\alpha+\beta)+ h^*(\gamma)
+[\alpha+\beta,\gamma]$ $$=h^*(\alpha+\beta)+
h^*(\alpha+\gamma)+h^*(\beta+\gamma)-h^*(\alpha)-h^*(\beta)-h^*(\gamma)
$$}

Finally it remains to show that $h^*$ respects the action of $G$.
Let $\alpha:S^2 \to X$ represent an element of $\pi_2(X)$.  We may
identify the complement of $+ \in S^2$ with the interior of the
closed unit disk $D^2 \subset \R^2$.  $\alpha$ extends naturally
to a map $\alpha':D^2 \to X$ which maps the boundary $\partial
D^2$ to the basepoint $* \in X$.  Let $R_\theta$ denote a positive
rotation of $D^2$ through the angle $\theta$.

We consider $h^*(\alpha):S^3 \to X$.  For each $(z,w) \in S^3$ we
have a fibre $\{(\lambda z, \lambda w) \vert\, \lambda
\bar{\lambda}=1\}$. If $z \neq 0$ this fibre completes a
revolution of the origin in the $w$- plane as it passes round the
$z$- plane.  The fibre $\{(0, \lambda)  \vert\, \lambda
\bar{\lambda}=1\}$ gets mapped to $* \in X$ by $h^*(\alpha)$.  We
may deform $h^*(\alpha)$ to thicken this fibre.  The remaining
fibres then form a solid torus.  We may therefore represent
$h^*(\alpha) \in \pi_3(X)$ by the map $\alpha^T:S^3 \to X$ defined
as follows:

\bigskip
Consider an unknotted solid torus $N \subset S^3$, not containing
the basepoint.  We make an untwisted topological identification of
$N$ with $S^1 \times D^2$.  Then:

\bigskip
$\alpha^T(\theta,d) =\alpha'(R_\theta(d)) \,\qquad \qquad
(\theta,d) \in N$

$\alpha^T(p)\quad=\,* \quad\quad\quad\qquad\qquad \qquad\, p\notin
N$

\bigskip
Let $I$ denote the closed unit interval and let $g:I \to X$
represent an element of $\pi_1(X)$, so $g(0)=g(1)=*$.  Then $g$
acts on $\alpha$ by stretching it round the loop $g$.  Formally we
have $\alpha g:D^2 \to X$ defined by:

\bigskip
$\alpha g(d)= \alpha'(2d) \qquad \qquad\qquad \quad \vert d \vert
\leq \frac{1}{2}$
\newline
${}$\quad \,\,\,$\alpha g(d)=g(2\vert d \vert -1) \,\quad \qquad
\qquad \vert d \vert \geq \frac{1}{2}$

Let $K \subset S^3$ be a thickened torus (parameterized $T^2
\times I$), bounding $N$ and not enclosing the basepoint.  Then
$h^*(\alpha g)$ is represented by the map $(\alpha g)^T:S^3 \to
X$, given by:

\bigskip
$(\alpha g)^T(\theta,d) =\alpha'(R_\theta(d)) \,\qquad \qquad
(\theta,d) \in N$

$(\alpha g)^T(t,\,i) \,= g(i) \,\qquad \qquad\,\,\,\qquad (t,\,i)
\,\in K$

$(\alpha g)^T(p)\quad=\,* \,\,\,\,\quad\qquad\qquad \qquad\,
p\,\notin N, K$

\bigskip
Let $L \subset S^3$ be a thickened sphere (parameterized $S^2
\times I$), enclosing $N$ and not enclosing the basepoint.  The
group element $g$ acts on an element of $\pi_3(X)$ by stretching
it round the loop $g$. Formally, $h^*(\alpha) g$ is represented by
the map $\alpha^Tg:S^3 \to X$, given by:

\bigskip
$\alpha^Tg(\theta,d) =\alpha'(R_\theta(d)) \,\qquad \qquad
(\theta,d) \in N$

$\alpha^Tg(s,\,i) \,= g(i) \,\qquad \qquad\,\,\,\qquad (s,\,i) \in
L$

$\alpha^Tg(p)\quad=\,* \,\,\,\,\quad\qquad\qquad \qquad\,
p\,\notin N, L$

\lem{Let $g$ be an element of $\pi_1(X)$ and let $\alpha$ be an
element of $\pi_2(X)$.  Then
$h^*(\alpha)g =h^*(\alpha g)$.}

\bigskip
Proof:  We must show that $\alpha^T g$ is homotopic to  $(\alpha
g)^T$.  Let $C \subset S^3$ be a solid cylinder which passes
through the hole in the solid torus $N$ and passes through the
thickened sphere $L$ on either side of this hole.

We may parameterize $C$ as $D^2 \times I$ so that for $(d,i) \in
C$, we have $$\alpha^T g (d,i) = m(i)$$ where $m=e \cdot g^{-1}
\cdot e \cdot g \cdot e$.  We define a homotopy $m_t$, $t \in I$,
from $m$ to $e$:

\bigskip
$m_t(i)= g(5r- .5) \qquad \qquad\quad r \in [.1, .3]$

$m_t(i)=* \,\,\,\,\,\qquad \qquad \qquad \qquad r \notin [.1, .3]$

\bigskip
\noindent where $r = \sqrt{t^2+(i- .5)^2}$.

\bigskip
We next construct a homotopy $f_t$ with $f_0=\alpha^T g$, which is
constant outside $C$.  Inside $C$ it is given by:

\bigskip

$f_t(d,i)=m_{(1-\vert d\vert)t}(i) $

\bigskip
Thus intuitively $f$ pulls together two regions of the thickened
sphere $L$ through the hole in the solid torus $N$, so that they
cancel out, leaving a thickened torus as required.

\hfill $\Box \,\,$

We may conclude:

\prop{$h^*$ is a $G$- quadratic map.}

Proof:  We get property i) from lemma 2.3, property ii) from lemma
2.6 and property iii) from lemma 2.7.

\hfill $\Box \,\,$

\sec{$\pi_3$ of a connected $2$- complex}

We now restrict to the case where $X$ is a connected 2- complex
with basepoint $*$ and $\pi_1(X)=G$.  The universal cover of $X$
is a simply connected 2- complex and hence a Moore space.

As $H_2(\tilde{X}; \Z)$ is a submodule of the free $\Z$- module
generated by the 2- cells of $\tilde{X}$, we have $H_2(\tilde{X};
\Z) \cong \Z^I$ for some set $I$ (see \cite{Hung}, Chap. IV,
Theorem 6.1).

From \cite{Baue}, Lemma 1.3.1, we know that up to homotopy there
is a unique CW complex of type $M(\Z^I,2)$.  Hence we have a
homotopy equivalence:
$$\bigvee_{i \in I} S^2_i \to \tilde{X}$$ preserving basepoint
(which we take to be $\bigcap_{i \in I} S^2_i$).

Composing with the covering map $\tilde{X} \to X$, we obtain a map
$m:\bigvee_{i} S^2_i \to X$. For $n>1$, the map $m$ induces an
isomorphism $\pi_n(\bigvee_{i } S^2_i) \to \pi_n(X)$.

\bigskip
\bigskip
For $j \in I$ let $w_j:S^2 \to \bigvee_{i } S^2_i$ be a
homeomorphism onto $S^2_j$, preserving basepoint.

\lem{i)  $\pi_2(\bigvee_{i } S^2_i)$ is freely generated over $\Z$
by the $w_i$. \newline ${}$\hspace{13mm} ii) $\pi_3(\bigvee_{i}
S^2_i)$ is freely generated over $\Z$ by the $h^*(w_i)$ and the
$[w_i, w_j]$, $i \neq j$.}

Proof: For i) see \cite{Hatc}, Example 4.26 and for ii) see
\cite{Hatc}, Example 4.52. \hfill $\Box$

\bigskip
\bigskip
Let $\alpha_i = m_*(w_i)$ for  $i \in I$.  As $m_*$ is an
isomorphism $\pi_2(\bigvee_{i } S^2_i) \to \pi_2(X)$, we know the
$\alpha_i$ freely generate $\pi_2(X)$ over $\Z$.

\lem{We have a $\Z$- basis for $\pi_3(X)$ given by the
$h^*(\alpha_i)$ and the $h^*(\alpha_i + \alpha_j)$, $i\neq j$.}

\proof{  As $m_*$ is an isomorphism $\pi_3(\bigvee_{i } S^2_i) \to
\pi_3(X)$, the $m_*(h^*(w_i))$ and the $m_*[w_i , w_j]$ freely
generate $\pi_3(X)$ over $\Z$. By lemma 2.2 and proposition 2.4i)
as well as lemma 2.5 we have:
\newline
\newline
$m_*(h^*(w_i)) =h^*(\alpha_i)$ \hspace{65mm}$\forall i$
\newline  $m_*[w_i, w_j]\,\,\, = [\alpha_i , \alpha_j] =
h^*(\alpha_i + \alpha_j) - h^*(\alpha_i)-h^*\alpha_j$
\hspace{16mm}$\forall i \neq j$}

\bigskip
\bigskip
By lemma 1.4 we have a $\ZG$- linear map $ \phi:S^2(\pi_2(X)) \to
\pi_3(X)$ which satisfies $\phi q =h^*$.

\lem{ $\phi$ is an isomorphism.}

\proof{Recall that the $q(\alpha_i)$ and the $q(\alpha_i +
\alpha_j)$, $i \neq j$ form a $\Z$- basis for $S^2(\pi_2(X))$. It
is sufficient to show that $\phi$ maps this basis to a $\Z$- basis
of $\pi_3(X)$.  By construction, $\phi q(\alpha_i) =
h^*(\alpha_i)$ and $\phi q(\alpha_i + \alpha_j) = h^*(\alpha_i +
\alpha_j)$.  By lemma 3.2 the $h^*(\alpha_i)$ and the
$h^*(\alpha_i + \alpha_j)$, $i \neq j$, form a $\Z$- basis of
$\pi_3(X)$ as required.}

From this lemma we may conclude our Theorem A:

\bigskip
\noindent{\bf Theorem A}  We have an isomorphism of $\ZG$ modules:
$\pi_3(X)\cong S^2(\pi_2(X))$

\sec{The effect of stabilizing $\pi_2$.}

  Let $X$ be a finite connected $2$- complex with finite fundamental
group $G$.  Let $X'$ be another finite connected $2$- complex with
fundamental group $G$.  By Schanuel's lemma we know that there
exist integers $a,b$ such that

$$
\pi_2(X) \oplus \ZG^a \cong \pi_2(X') \oplus \ZG^b
$$

  In this section we investigate the corresponding relationship
between $\pi_3(X)$ and $\pi_3(X')$.

  If we let $J= \pi_2(X)$ and $J'= \pi_2(X')$, then from Theorem A we have
\newline
$\pi_3(X) \cong S^2(J)$ and $\pi_3(X')\cong S^2(J')$.  We know
that

$$
S^2(J\oplus \ZG^a) \cong S^2(J' \oplus \ZG^b)
$$

  The next few lemmas give us an expansion of this.

\lem{Let $A_i$, $i=1, \cdots, n$, be modules over $\ZG$ with
finitely generated free underlying Abelian groups.  Then
$$
S^2(\bigoplus_{i=1}^n A_i)= \bigoplus_{i=1}^n S^2(A_i) \oplus \bigoplus_{i <j} A_i
\otimes_\Z A_j
$$
}

Proof:  {For each $i$, let the $\vec{e_{i,r}}$ be a basis over
$\Z$ for $A_i$.  Then $S^2(\bigoplus_{i=1}^n A_i)$ is freely
generated over $\Z$, by elements of the form:
\newline
\newline
${}$ \hspace{7mm}$\vec{e_{i,r}} \otimes \vec{e_{i,r}}$,
\newline
${}$ \hspace{7mm}$\vec{e_{i,r}} \otimes
\vec{e_{i,s}}+\vec{e_{i,s}} \otimes \vec{e_{i,r}}\qquad\qquad
\qquad r<s$,
\newline
${}$ \hspace{7mm}$\vec{e_{i,r}} \otimes
\vec{e_{j,s}}+\vec{e_{j,s}} \otimes \vec{e_{i,r}}\qquad \qquad
\qquad i<j$,
\newline
\newline
${}$ \hspace{2mm}  For each $i$, the $\Z$- linear span of the
\newline
\newline
${}$ \hspace{7mm}$\vec{e_{i,r}} \otimes \vec{e_{i,r}}$,
\newline
${}$ \hspace{7mm}$\vec{e_{i,r}} \otimes
\vec{e_{i,s}}+\vec{e_{i,s}} \otimes \vec{e_{i,r}}\qquad\qquad
\qquad r<s$,
\newline
\newline
is closed under the group action and is isomorphic (over $\ZG$) to
$S^2(A_i)$.
\newline
\newline
${}$\hspace{2mm}Similarly, for each pair $i,j$, with $i <j$, the
$\Z$- linear span of the
\newline
\newline
${}$ \hspace{7mm}$\vec{e_{i,r}} \otimes
\vec{e_{j,s}}+\vec{e_{j,s}} \otimes \vec{e_{i,r}}\qquad \qquad
\qquad \forall r,s$
\newline
\newline
is closed under the group action. Furthermore we have an
isomorphism (over $\ZG$) from it to $A_i \otimes_\Z A_j$, which
maps:
$$\vec{e_{i,r}} \otimes \vec{e_{j,s}}+\vec{e_{j,s}} \otimes \vec{e_{i,r}} \quad\mapsto\quad \vec{e_{i,r}} \otimes
\vec{e_{j,s}}$$ ${}$ \hspace{2mm}So over $\ZG$, the module
$S^2(\bigoplus_{i=1}^n A_i)$ decomposes into the sum:
$$
\bigoplus_{i=1}^n S^2(A_i) \oplus \bigoplus_{i <j} A_i
\otimes_\Z A_j
$$\hfill $\Box \,\,$}

\lem{Let $A$ be a $\ZG$ module whose underlying Abelian group is isomorphic to
$\Z^k$.  Then $A \otimes_\Z \ZG \cong \ZG^k$.}

Proof:  The action of an element $g \in G$ on $A$ is a $\Z$-
linear isomorphism.  Hence, if $\vec{e_1}, \cdots, \vec{e_k}$ is a
$\Z$- linear basis for $A$, then so is $\vec{e_1}g, \cdots,
\vec{e_k}g$. Hence $A \otimes_\Z \ZG$ is freely generated over
$\Z$ by elements of the form $\vec{e_i}g \otimes g$, for $g\in G$
and $i \in \{1,\cdots,k\}$.

  For each $i$, the $\Z$ linear span
of the $\vec{e_i}g \otimes g$, $g\in G$, is closed under the
action of $G$. Furthermore, we have a $\ZG$- linear isomorphism
from it to $\ZG$, which maps $\vec{e_i}g \otimes g \mapsto g$.

Hence,
$$
A \otimes_\Z \ZG \cong \bigoplus_{i=1}^k\ZG \cong \ZG^k
$$

\hfill $\Box \,\,$

\bigskip
Let $p$ denote the number of pairs, $\{g, g^{-1}\} \in G$, with $g
\neq g^{-1}$.  Let $T$ denote the set of $g \in G$ with $g$ having order $2$.

\define{} For a finite group $G$, we set

$$V_G = \bigoplus_{t \in T} (1+t) \ZG
$$

\lem{$S^2(\ZG) \cong \ZG^{1+p} \oplus V_G$}

Proof:  Let the the elements of $G$ be ordered.  We have a $\Z$-
basis for $S^2(\ZG)$ given by:
\newline
\newline
${}$ \hspace{7mm}$(e \otimes e)h$ \hfill  $h \in G \qquad \qquad
\qquad\qquad \qquad$
\newline
${}$ \hspace{7mm}$(e \otimes g + g\otimes e)h$\hfill  $g < g^{-1},
\, \quad h \in G\qquad \qquad \qquad\qquad\qquad$
\newline
${}$ \hspace{7mm}$(e \otimes t + t\otimes e)h$ \hfill  $t \in
T,\quad  h < th\qquad \qquad \qquad\qquad\qquad$
\newline
\newline
${}$\hspace{2mm} The $\Z$- linear span of the $\{(e \otimes e)h
\vert \,\, h \in G\}$ is isomorphic (over $\ZG$) to $\ZG$.  Also,
for $g < g^{-1}$, the $\Z$- linear span of the $\{(e \otimes g +
g\otimes e)h\vert \,\, h \in G\}$ is isomorphic (over $\ZG$) to
$\ZG$.

Given $t \in T$, the $\Z$- linear span of the $\{(e \otimes t +
t\otimes e)h\vert \,\, h < th \}$ is isomorphic to $(e \otimes
t)(1+t) \ZG \subset \ZG \otimes_\Z \ZG$.

Finally note that as $(e \otimes t)\ZG \cong \ZG)$, we have $(e
\otimes t)(1+t) \ZG \cong (1+t)\ZG$.

\hfill $\Box \,\,$

\bigskip
 Note that for $\ZG$- modules $A$, $B$, we have a $\ZG$- linear
isomorphism
\newline
$A \otimes_\Z B \to B \otimes_\Z A$, given by sending $\vec{a}
\otimes \vec{b} \mapsto \vec{b} \otimes \vec{a}$.

\bigskip
Using lemma 4.1 to expand the identity $S^2(J\oplus \ZG^a) \cong
S^2(J' \oplus \ZG^b)$ we obtain:

$$
\quad \, \,S^2(J) \oplus S^2(\ZG)^a \oplus (J \otimes_\Z \ZG)^a
\oplus (\ZG \otimes_\Z \ZG)^{a(a-1)/2}
$$
$$
\cong S^2(J') \oplus S^2(\ZG)^b \oplus (J' \otimes_\Z \ZG)^b
\oplus (\ZG \otimes_\Z \ZG)^{b(b-1)/2}
$$

\bigskip Let $k$ denote the $\Z$- rank of $J$, and let $k'$ denote
the $\Z$- rank of $J'$.  Also let $n$ denote the order of $G$.  By
Theorem A we have $S^2(J) \cong \pi_3(X)$ and $S^2(J') \cong
\pi_3(X')$. Applying lemmas 4.2 and 4.4 to the equation above
gives:

$$
\quad\,\pi_3(X) \oplus (\ZG^{1+p} \oplus V_G)^a \oplus \ZG^{ka}
\oplus \ZG^{na(a-1)/2}
$$
$$
\cong \pi_3(X') \oplus (\ZG^{1+p} \oplus V_G)^b \oplus \ZG^{k'b}
\oplus \ZG^{nb(b-1)/2}
$$

So:
$$
\pi_3(X) \oplus \ZG^{(1+p+k+n(a-1)/2)a} \oplus {V_G}^a \cong
\pi_3(X') \oplus \ZG^{(1+p+k'+n(b-1)/2)b} \oplus {V_G}^b
$$

\thm{Let $X$ and $X'$ be finite connected $2$- complexes, with
finite fundamental group $G$. Then there exist integers $a,b,q,r$
such that:
$$
\pi_3(X) \oplus \ZG^q \oplus {V_G}^a  \cong \pi_3(X') \oplus \ZG^r
\oplus {V_G}^b
$$
}

\bigskip
Note that if the order of $G$ is odd then it does not contain any
elements of order $2$.  Hence $V_G=0$ and we have obtained Theorem
B:

\bigskip
\noindent {\bf Theorem B} {\it Let $X$ and $X'$ be finite
connected 2- complexes, with finite fundamental group $G$.  If the
order of $G$ is odd, then $\pi_3(X)$ and $\pi_3(X')$ are stably
equivalent.}

\sec{Calculating $S^2(IG^*$)}

Again let $X$ be a finite connected $2$- complex with finite
fundamental group $G$.  Rationally we have $\pi_2(X) \otimes \Q
\cong IG^* \otimes \Q \oplus \QG^m$, for some $m \geq 0$.
Consequently $\pi_3(X) \otimes \Q \cong S^2(IG^* \oplus \ZG^m)
\otimes \Q$ as modules over $\QG$.  For this reason we examine
$S^2(IG^*)$.

We will then conclude by considering some examples where
$\pi_2(X)$ actually equals $IG^*$.  For this to happen the group
$G$ must necessarily have a balanced presentation and free period
4.

\bigskip
$IG^*$ is defined by the short exact sequence:

$$
0 \to \Z \to \ZG \stackrel{\delta}{\to} IG^*\to0
$$
where $1 \in \Z$ is mapped to $\sum_{g \in G} g \in \ZG$, denoted
$\Sigma$.

\bigskip
We now consider the $\ZG$- linear surjection $\delta':S^2(\ZG) \to
S^2(IG^*)$ defined by: $\,\,\qquad\qquad\qquad \qquad g \otimes g
\qquad\, \mapsto \quad \,\,\delta(g) \otimes \delta(g)$

$\quad \qquad \qquad  \qquad\qquad \qquad g \otimes h +h \otimes g
\mapsto \delta(g) \otimes \delta(h) + \delta(h) \otimes \delta(g)$

\lem{The kernel of $\delta'$ is generated over $\ZG$, by $\Sigma
\otimes e + e \otimes \Sigma$ and $\Sigma \otimes \Sigma$.}

Proof:  We may take a $\Z$- basis for $\ZG$ given by $\{g \in G
\vert \,\, g \neq e\} \sqcup \{\Sigma\}$.  We have a $\Z$- basis
for $IG^*$ given by $\{p(g)  \vert \,\, g \in G, g \neq e \}$.

These basis' induce basis' on $S^2(\ZG)$ and $S^2(IG^*)$
respectively. We have:

\bigskip
$\delta'(\Sigma \otimes \Sigma)=0$

$\delta'(g \otimes \Sigma + \Sigma \otimes g)=0 \qquad \qquad g
\neq e$

\bigskip
\noindent and $\delta'$ maps the other elements of the induced
basis for $S^2(\ZG)$ to elements of the induced basis for
$S^2(IG^*)$.

Hence the kernel of $\delta'$ is generated over $\Z$ by $\Sigma
\otimes \Sigma$ and $g \otimes \Sigma + \Sigma \otimes g$, $g \neq
e$. Clearly this is contained in the $\ZG$ linear span of $\Sigma
\otimes e + e \otimes \Sigma$ and $\Sigma \otimes \Sigma$.

\hfill $\Box \,\,$

\bigskip
Let $u = e \otimes \Sigma + \Sigma \otimes e$.  Then we have the
kernel of $\delta'$ generated over $\ZG$ by $\{u, u\Sigma/2\}$.
Note that $u\ZG \cong \ZG$.

As $u\Sigma/2$ is in the $\QG$- linear span of $u$, the kernel of
$\delta' \otimes \Q:S^2(\QG) \to S^2(IG^*) \otimes \Q$ is
generated by $u$.  We therefore have a short exact sequence:
$$
0 \to \QG \to S^2(\QG) \to S^2(IG^*) \otimes \Q \to 0
$$

As surjections over $\QG$ split and cancellation of finitely
generated modules holds over $\QG$, we may write:
$$
 S^2(IG^*) \otimes \Q \cong  S^2(\QG) / \QG
$$

\bigskip
Rationally, for some $a \in \Z$, $\pi_2(X) \otimes \Q \cong (IG^*
\otimes \Q) \oplus \QG^a$.  Consequentially, we have:
$$
\pi_3(X) \otimes \Q  \cong S^2(\pi_2(X)) \otimes \Q \cong
S^2(\pi_2(X) \otimes \Q) \cong S^2(IG^* \oplus \ZG^a) \otimes \Q$$

\bigskip
From lemmas 4.1, 4.2 and 4.4, for some $q \in \Z$ we have:
$$ S^2(IG^* \oplus \ZG^a) \otimes \Q \cong (S^2(IG^*) \otimes Q) \oplus
\QG^q \oplus (V_G \otimes \Q)^a$$

As $S^2(IG^*)\otimes \Q \cong (S^2(\ZG) \otimes \Q) / \QG$, from
lemma 4.4 we may conclude:

\thm{
$$
\pi_3(X) \otimes \Q \cong \QG^{p+q} \oplus (V_G \otimes \Q)^{a+1}
$$
}

Again if the order of $G$ is odd then $V_G=0$ so we obtain our
Theorem C:

\bigskip
\noindent {\bf Theorem C} {\it Let $X$ be a finite connected $2$-
complex with finite fundamental group $G$.  If the order of $G$ is
odd, then $\pi_3(X) \otimes \Q$ is free over $\QG$.}

\bigskip
 Let $S$ be a subset of $G$, containing precisely one of $g$ or
$g^{-1}$, for each pair $g, g^{-1}$, with $g \neq g^{-1}$.  Let
$\vec{e_g} = g \otimes e + e \otimes g$ for each $g \in S$.  Let
$\vec{e_t}=t \otimes e + e \otimes t$ for each $t \in T$.  Let
$\vec{e_0}= e \otimes e$.  Then from the proof of lemma 4.4 we
have:
$$
S^2(\ZG)= \vec{e_0}\ZG \oplus \bigoplus_{g \in S} {e_g}\ZG \oplus
\bigoplus_{t \in T} \vec{e_t}\ZG
$$
Also,
$$
u = e\otimes \Sigma + \Sigma \otimes e = \vec{e_0} 2+\sum_{g \in
S}\vec{e_g}(1+g^{-1}) + \sum_{t\in T}\vec{e_t}
$$

 \noindent and
$$
S^2(IG^*)= \frac {S^2(\ZG) } { \quad <u, u { \frac{\Sigma}{2}}>}
$$

The relation $u=0$ is equivalent to $\vec{e_0} 2=-(\sum_{g \in
S}\vec{e_g}(1+g^{-1}) + \sum_{t\in T}\vec{e_t})$.

\bigskip
For $t \in T$, we let $S_t \in \ZG$ be some element satisfying
$(1+t)S_t = \Sigma$.  Hence we have $\vec{e_t} S_t = (e \otimes t)
(1+t)S_t = (e \otimes t) \Sigma = (e \otimes t) (1+t)\Sigma /2 =
\vec{e_t} \Sigma /2$.

\bigskip
The relation $u\Sigma/2=0$ is equivalent to $\vec{e_0} \Sigma
=-(\sum_{g \in S}\vec{e_g}\Sigma + \sum_{t\in T}\vec{e_t} S_t)$.

\bigskip
\noindent Let
$$M=\bigoplus_{g \in S} \vec{e_g}\ZG \oplus
\bigoplus_{t \in T} \vec{e_t}\ZG
$$

\bigskip
\noindent and let $\vec{u_M}  = -(\sum_{g \in
S}\vec{e_g}(1+g^{-1}) + \sum_{t\in T}\vec{e_t})$.

\prop{If $\pi_2{X}=IG^*$ then $\pi_3(X)=M[\frac{\vec{u_M}}{2}]$.}

\proof{We need to check that $\frac{\vec{u_M}}{2}$ satisfies the
relations which $\vec{e_0}$ was subject to:
\newline
\newline ${}$\hspace{10mm}
$\frac{\vec{u_M}}{2}2 \,= -(\sum_{g \in S}\vec{e_g}(1+g^{-1}) +
\sum_{t\in T}\vec{e_t})$
\newline
\newline ${}$\hspace{10mm}
$\frac{\vec{u_M}}{2} \Sigma =-(\sum_{g \in S}\vec{e_g}\Sigma +
\sum_{t\in T}\vec{e_t} S_t)$}

\bigskip
\egg{$G=C_3$ and $X$ is the Cayley complex associated to the
balanced presentation: $<\!x \vert\,\, x^3=e \!>$.}

From \cite{John}, 41(i) we have $\pi_2(X)=IG^*$.  Hence:
$$\pi_3(X)  \cong S^2(IG^*) \cong \vec{e_{x^2}}\ZG[(1+x)/2] \cong \ZG[\frac{(1+x)}{2}]
$$

\bigskip
\bigskip
\egg{$G=Q_8$ and $X$ is the Cayley complex associated to the
balanced presentation: $<\!x,\, y \vert\,\, x^2= y^{2},\, xyx=y
\!>$.}

Again from \cite{John}, 41(iii) we have $\pi_2(X)=IG^*$.  There is
only one element of order two in $Q_8$ so we have:
$$
\pi_3(X) \cong S^2(IG^*) \cong (\vec{e_x} \ZG \oplus \vec{e_y} \ZG
\oplus \vec{e_{xy}} \ZG \oplus \vec{e_{y^2}}
\ZG)[\frac{\vec{u_M}}{2}]
$$
where
$\vec{u_M}=-(\vec{e_x}(1+x^3)+\vec{e_y}(1+y^3)+\vec{e_{xy}}(1+xy^3)+\vec{e_{y^2}})$.

Recall from the proof of lemma 4.4 that:  $$\vec{e_x} \ZG \cong
\vec{e_y} \ZG \cong \vec{e_{xy}} \ZG \cong \ZG, \qquad \qquad
\vec{e_{y^2}}\ZG \cong (1+y^2)\ZG
$$


\begin{thebibliography}{}

\bibitem{Baue} Hans-Joachim Baues, {\sl Homotopy Type and Homology} (OUP, 1996)

\bibitem{Hatc} Allen Hatcher, {\sl Algebraic Topology} (CUP, 2002)

\bibitem{Hung} Thomas W. Hungerford, {\sl Algebra} (GTM 73, 1974)

\bibitem{John} F.E.A. Johnson,   {\sl Stable Modules and the D(2)- Problem}  (LMS 301, 2003)

\bibitem{Whit} George W. Whitehead, {\sl Elements of Homotopy Theory} (GTM 61, 1978)





%\bibitem{Kirb} R Kirby, L Siebenmann ; {\sl Foundational essays on topological manifolds, smoothings and triangulations} : Ann. Math. Studies
% 88, Princeton Univ. Press, Princeton, NJ (1977)

%\bibitem{Macl} Saunders Mac Lane ; {\sl Homology} : Springer-Verlag (1963)

%\bibitem{May} J.P.May ; {\sl A Concise Course in Algebraic Topology} : Chicago
%Lectures in Mathematics (1999)

%\bibitem{Miln1} John Milnor ; {\sl On the construction FK} : Algebraic Topology - a student's guide (J F Adams) :
%LMS Lecture Note Series no. 4, CUP (1972)

%\bibitem{Vick} Vick, James W.: {\sl Homology Theory} : GTM 145 (1994)

\end{thebibliography}
\end{document}